
%
%
%
%
%
%
\magnification=\magstephalf      
%
%
\vsize=7.5truein                 
\hsize=5.2truein                 
\newskip\stdskip                 
\stdskip=6pt plus3pt minus3pt    
\medskipamount=\stdskip          
\parindent=0pt                   
\parskip=\stdskip                
\abovedisplayskip=\stdskip       
\belowdisplayskip=\stdskip       
\mathsurround=0.75pt             
\overfullrule=0pt                
%
%
\def\ppar{\par\goodbreak\vskip 8pt plus 4pt minus 4pt}     
%
%
\def\stdspace{\hskip 0.75em plus 0.15em\ignorespaces}
\let\qua\stdspace 
%
%
%
%
%
%
%
\def\hexnumber#1{\ifcase#1 0\or 1\or 2\or 3\or 4\or 5\or 6\or 7\or 8\or
 9\or A\or B\or C\or D\or E\or F\fi}
%
%
\font\thirtnmsa=msam10 scaled 1315    
\font\tenmsa=msam10          \font\ninemsa=msam9
\font\sevenmsa=msam7         \font\sixmsa=msam6
\font\fivemsa=msam5
%
%
\newfam\msafam                  \textfont\msafam=\tenmsa
\scriptfont\msafam=\sevenmsa    \scriptscriptfont\msafam=\fivemsa
\edef\hexa{\hexnumber\msafam}        
\def\msa{\fam\msafam\tenmsa}         
%
%
\font\thirtnmsb=msbm10 scaled 1315   
\font\tenmsb=msbm10      \font\ninemsb=msbm9
\font\sevenmsb=msbm7     \font\sixmsb=msbm6
\font\fivemsb=msbm5
%
\newfam\msbfam                   \textfont\msbfam=\tenmsb       
\scriptfont\msbfam=\sevenmsb     \scriptscriptfont\msbfam=\fivemsb
\edef\hexb{\hexnumber\msbfam}    
\def\msb{\fam\msbfam\tenmsb}     
%
%
\font\thirtneufm=eufm10 scaled 1315   
\font\teneufm=eufm10                 \font\nineeufm=eufm9
\font\seveneufm=eufm7                \font\sixeufm=eufm6
\font\fiveeufm=eufm5
%
\newfam\eufmfam                    \textfont\eufmfam=\teneufm
\scriptfont\eufmfam=\seveneufm     \scriptscriptfont\eufmfam=\fiveeufm
\edef\hexf{\hexnumber\eufmfam}      
\def\frak{\fam\eufmfam\teneufm}     
%
%
%
\font\thirtnrm=cmr10 scaled 1315    
\font\ninerm=cmr9                   \font\sixrm=cmr6   
%
\font\thirtni=cmmi10 scaled 1315    
\font\ninei=cmmi9                   \font\sixi=cmmi6  
%
\font\thirtnsy=cmsy10 scaled 1315   
\font\ninesy=cmsy9                  \font\sixsy=cmsy6  
%
\font\thirtnbf=cmbx10 scaled 1315   
\font\ninebf=cmbx9                  \font\sixbf=cmbx6  
%
%
\font\thirtnex=cmex10 scaled 1315   
\font\nineex=cmex9                  
%
%
\font\thirtnit=cmti10 scaled 1315  
\font\nineit=cmti9                  
%
\font\thirtnsl=cmsl10 scaled 1315  
\font\ninesl=cmsl9                  
%
\font\thirtntt=cmtt10 scaled 1315  
\font\ninett=cmtt9                  
%
%
%
%
\def\small{%
%
%
\textfont0=\ninerm \scriptfont0=\sixrm \scriptscriptfont0=\fiverm
\def\rm{\fam0\ninerm}
%
%
\textfont1=\ninei \scriptfont1=\sixi \scriptscriptfont1=\fivei
%
%
\textfont2=\ninesy \scriptfont2=\sixsy \scriptscriptfont2=\fivesy
%
%
\textfont3=\nineex \scriptfont3=\nineex \scriptscriptfont3=\nineex
%
%
\textfont\bffam=\ninebf \scriptfont\bffam=\sixbf
\scriptscriptfont\bffam=\fivebf \def\bf{\fam\bffam\ninebf}%
%
%
\textfont\itfam=\nineit \def\it{\fam\itfam\nineit}%
\textfont\slfam=\ninesl \def\sl{\fam\slfam\ninesl}%
\textfont\ttfam=\ninett \def\tt{\fam\ttfam\ninett}%
%
%
%
\textfont\msafam=\ninemsa \scriptfont\msafam=\sixmsa
\scriptscriptfont\msafam=\fivemsa \def\msa{\fam\msafam\ninemsa}%
%
%
\textfont\msbfam=\ninemsb \scriptfont\msbfam=\sixmsb
\scriptscriptfont\msbfam=\fivemsb \def\msb{\fam\msbfam\ninemsb}%
%
%
\textfont\eufmfam=\nineeufm  \scriptfont\eufmfam=\sixeufm
\scriptscriptfont\eufmfam=\fiveeufm \def\frak{\fam\eufmfam\nineeufm}%
%
%
%
\normalbaselineskip=11pt%
\setbox\strutbox=\hbox{\vrule height8pt depth3pt width0pt}%
%
%
\normalbaselines\rm
%
%
\stdskip=4pt plus2pt minus2pt    
\medskipamount=\stdskip          
\parskip=\stdskip                
\abovedisplayskip=\stdskip       
\belowdisplayskip=\stdskip       
\def\ppar{\par\goodbreak\vskip 6pt plus 3pt minus 3pt}%
%
%
\def\section##1{\global\advance\sectionnumber by 1
\vskip-\lastskip\penalty-800\vskip 20pt plus10pt minus5pt 
\egroup{\bf\number\sectionnumber\quad##1}\bgroup\small         
\vskip 6pt plus3pt minus3pt
\nobreak\resultnumber=1}
}    
%
\def\beginsmall{\bgroup\small}
\let\endsmall\egroup
%
%
%
%
\def\large{%
\textfont0=\thirtnrm \scriptfont0=\ninerm \scriptscriptfont0=\sevenrm
\def\rm{\fam0\thirtnrm}%
\textfont1=\thirtni \scriptfont1=\ninei \scriptscriptfont1=\seveni
\textfont2=\thirtnsy \scriptfont2=\ninesy \scriptscriptfont2=\sevensy
\textfont3=\thirtnex \scriptfont3=\thirtnex \scriptscriptfont3=\thirtnex
\textfont\bffam=\thirtnbf \scriptfont\bffam=\ninebf
\scriptscriptfont\bffam=\sevenbf \def\bf{\fam\bffam\thirtnbf}%
\textfont\itfam=\thirtnit \def\it{\fam\itfam\thirtnit}%
\textfont\slfam=\thirtnsl \def\sl{\fam\slfam\thirtnsl}%
\textfont\ttfam=\thirtntt \def\tt{\fam\ttfam\thirtntt}%
\textfont\msafam=\thirtnmsa \scriptfont\msafam=\ninemsa
\scriptscriptfont\msafam=\sevenmsa \def\msa{\fam\msafam\thirtnmsa}%
\textfont\msbfam=\thirtnmsb \scriptfont\msbfam=\ninemsb
\scriptscriptfont\msbfam=\sevenmsb \def\msb{\fam\msbfam\thirtnmsb}%
\textfont\eufmfam=\thirtneufm  \scriptfont\eufmfam=\nineeufm
\scriptscriptfont\eufmfam=\seveneufm \def\frak{\fam\eufmfam\teneufm}%
\normalbaselineskip=16pt%
\setbox\strutbox=\hbox{\vrule height11.5pt depth4.5pt width0pt}%
\normalbaselines\rm}%
\let\Large\large   
%

%

%
\mathchardef\plussquare="0\hexa01
\mathchardef\nge="3\hexb0B
\mathchardef\maltesecross="0\hexa7A
\mathchardef\del="0\hexf01
%
%
%
%
\font\sc=cmcsc10
%
%
%
%
\def\sqr#1#2{{\vcenter{\vbox{\hrule  height.#2truept
	\hbox{\vrule width.#2truept height#1truept 
	\kern#1truept \vrule width.#2truept}
	\hrule height.#2truept}}}}
\def\sq{\sqr55}    
%
%
%
%
\newcount\sectionnumber            
\newcount\resultnumber             
\sectionnumber=0\resultnumber=1    
%
%
%
\def\section#1{\global\advance\sectionnumber by 1
\xdef\nextkey{\number\sectionnumber}
\vskip-\lastskip\penalty-800\vskip 20pt plus10pt minus5pt 
{\large\bf\number\sectionnumber\quad#1}         
\vskip 8pt plus4pt minus4pt
\nobreak\resultnumber=1}      
%
%
%
%
%
         
%
%
%
%

%
\def\proc#1{\xdef\nextkey{\number\sectionnumber.\number\resultnumber}%
\vskip-\lastskip\ppar\bf%
\noindent#1\ \number\sectionnumber.\number\resultnumber
\stdspace\sl\global\advance\resultnumber by 1\ignorespaces}
 
%
%
\def\proof#1{\vskip-\lastskip\ppar\noindent{\bf#1}%
\stdspace\rm\ignorespaces}        
\let\endproof\endprf              
%
%
%
%
%
%
%
%
\def\proclaim#1{\vskip-\lastskip\ppar\bf%
\noindent#1\stdspace\sl\ignorespaces} 
\let\endproclaim\endproc
%
%
%
%
\def\rk#1{\vskip-\lastskip\ppar{\bf #1}\stdspace\ignorespaces}                

%
%
%
%
%
%
\def\label{\xdef\nextkey{\number\sectionnumber.\number\resultnumber}%
\number\sectionnumber.\number\resultnumber
\global\advance\resultnumber by 1}
%
%
%
%
%
%
%
%
%
%
%
%
%
%
%
%
\newcount\refnumber              
\refnumber=1                     
\long\def\reflist#1\endreflist{%
\long\def\thereflist{#1}{\def\refkey##1##2\par{\xdef##1{\number\refnumber}%
\global\advance\refnumber by 1}%
\def\key##1##2\par{\expandafter\xdef%
\csname##1\endcsname{\number\refnumber}%
\global\advance\refnumber by 1}#1\par}}
\long\def\references{%
\penalty-800\vskip-\lastskip\vskip 15pt plus10pt minus5pt 
{\large\bf References}\ppar 
{\leftskip=25pt\frenchspacing    
\small\parskip=3pt plus2pt       
\def\refkey##1##2\par{\noindent  
\llap{[##1]\stdspace}\ignorespaces##2\par}         
\def\key##1##2\par{\noindent  
\llap{[\ref{##1}]\stdspace}\ignorespaces##2\par}  
\def\,{\thinspace}\thereflist\par}}
%
%
%
\newcount\footnotenumber         
\footnotenumber=1                
\def\fnote#1{\xdef\nextkey{\number\footnotenumber}%
{\small\ifnum\footnotenumber>9\parindent=14pt%
\else\parindent=10pt\fi\footnote{$^{\number\footnotenumber}$}%
{\hglue-5pt#1}\global\advance\footnotenumber by 1}}
%
%
%
%
%
%
%
\newcount\figurenumber          
\figurenumber=1                 
\def\caption#1{\xdef\nextkey{\number\figurenumber}%
\cl{\small Figure \number\figurenumber: #1}%
\global\advance\figurenumber by 1}
\def\figurelabel{\xdef\nextkey{\number\figurenumber}%
\cl{\small Figure \number\figurenumber}%
\global\advance\figurenumber by 1}
\long\def\figure#1\endfigure{{\xdef\nextkey{\number\figurenumber}%
\let\captiontext\relax\def\caption##1{\xdef\captiontext{##1}}%
\midinsert\cl{\ignorespaces#1\unskip\unskip\unskip\unskip}\vglue6pt\cl{\small 
Figure \number\figurenumber\ifx\captiontext\relax\else: \captiontext
\fi}\endinsert\global\advance\figurenumber by 1}}
%
%
%
%
%
%
%
\def\nextkey{??}   
%
\def\key#1{\expandafter\xdef\csname #1\endcsname{\nextkey}}
\def\ref#1{\expandafter\ifx\csname #1\endcsname\relax
\immediate\write16{Reference {#1} undefined}??\else
\csname #1\endcsname\fi}
%
%
%
%
%
%
%
\newread\gtinfile
\newwrite\gtreffile
\def\useforwardrefs{
\openin\gtinfile\jobname.ref
\ifeof\gtinfile
\closein\gtinfile
\immediate\write16{No file \jobname.ref}
\else
\closein\gtinfile
\input \jobname.ref
\fi
\immediate\openout\gtreffile \jobname.ref
%
%
\def\key##1{{\def\\{\noexpand}%
\expandafter\xdef\csname ##1\endcsname{\nextkey}%
\immediate\write\gtreffile{\\\expandafter\\\def\\\csname ##1\\\endcsname%
{\nextkey}}}}
%
%
\long\def\reflist##1\endreflist{%
\long\def\thereflist{##1}{\def\refkey####1####2\par{\xdef####1{%
\number\refnumber}{\def\\{\noexpand}\immediate\write\gtreffile
{\\\def\\####1{\number\refnumber}}}\global\advance\refnumber by 1}%
\def\key####1####2\par{\expandafter\xdef%
\csname####1\endcsname{\number\refnumber}%
{\def\\{\noexpand}\immediate\write\gtreffile
{\\\expandafter\\\def\\\csname ####1\\\endcsname{\number\refnumber}}}
\global\advance\refnumber by 1}##1\par}}
\long\def\biblio##1\endbiblio{\reflist##1\endreflist\references}%
%
%
\def\numkey##1{{\def\\{\noexpand}%
\xdef##1{\number\sectionnumber.\number\resultnumber}
\immediate\write\gtreffile{\\\def\\##1%
{\number\sectionnumber.\number\resultnumber}}}}
\def\seckey##1{{\def\\{\noexpand}\xdef##1{\number\sectionnumber}
\immediate\write\gtreffile{\\\def\\##1{\number\sectionnumber}}}}
\def\figkey##1{\xdef##1{\number\figurenumber}%
{\def\\{\noexpand}\immediate\write\gtreffile%
{\\\def\\##1{\number\figurenumber}}}
\number\figurenumber\global\advance\figurenumber by 1}
}   
%
%
%
%
\def\figkey#1{\xdef#1{\number\figurenumber}%
\number\figurenumber\global\advance\figurenumber by 1}
\def\fig#1#2\endfig{%
\midinsert\cl{#2}\vglue6pt\cl{\small Figure #1}\endinsert}
\def\newfig{\number\figurenumber\global\advance\figurenumber by 1}
\def\numkey#1{\xdef#1{\number\sectionnumber.\number\resultnumber}}
\def\seckey#1{\xdef#1{\number\sectionnumber}}
%
%
%
%
%
%
%
%
%
\def\verb{\catcode`\"=\active}       
\def\brev{\catcode`\"=12}            
\brev                                
\verb                                
{\obeyspaces\gdef {\ }}              
{\catcode`\`=\active\gdef`{\relax\lq}}
\def"{%
\begingroup\baselineskip=12pt\def\par{\leavevmode\endgraf}%
\tt\obeylines\obeyspaces\parskip=0pt\parindent=0pt%
\catcode`\$=12\catcode`\&=12\catcode`\^=12\catcode`\#=12%
\catcode`\_=12\catcode`\~=12%
\catcode`\{=12\catcode`\}=12\catcode`\%=12\catcode`\\=12%
\catcode`\`=\active\let"\endgroup}
\brev      
%
%
%
%
%
%
\def\item#1{\par\leavevmode\llap{#1\stdspace}%
\ignorespaces}                             
%
%

%
%
\def\np{\vfil\eject}         
\def\nl{\hfil\break}         
\def\cl{\centerline}         
\def\agt{{\mathsurround=0pt\it$\cal A\mskip-.7mu$lgebraic \&\ 
$\cal G\mskip-2mu$eometric $\cal T\!\!$opology}}  
%
%
%

%
%
%
%
%
\def\title#1{\def\thetitle{#1}}
\def\shorttitle#1{\def\theshorttitle{#1}}
\def\author#1{\edef\previousauthors{\theauthors}
 \ifx\theauthors\relax\def\theauthors{#1}\else
 \def\theauthors{\previousauthors\par#1}\fi}

\let\authors\author        
\def\address#1{\edef\previousaddresses{\theaddress}
 \ifx\theaddress\relax\def\theaddress{#1}\else
 \def\theaddress{\previousaddresses\par\vskip 2pt\par#1}\fi}
\def\secondaddress#1{\edef\previousaddresses{\theaddress}
 \ifx\theaddress\relax\def\theaddress{#1}\else
 \def\theaddress{\previousaddresses\par{\rm and}\par#1}\fi}   

\def\email#1{\edef\previousemails{\theemail}
 \ifx\theemail\relax\def\theemail{#1}\else
 \def\theemail{\previousemails\hskip 0.75em\relax#1}\fi}
\def\secondemail#1{\edef\previousemails{\theemail}
 \ifx\theemail\relax\def\theemail{#1}\else
 \def\theemail{\previousemails\hskip 0.75em{\rm and}\hskip 0.75em
 \relax#1}\fi}
\def\url#1{\edef\previousurls{\theurl}
 \ifx\theurl\relax\def\theurl{#1}\else
 \def\theurl{\previousurls\hskip 0.75em\relax#1}\fi}
\def\secondurl#1{\edef\previousurls{\theurl}
 \ifx\theurl\relax\def\theurl{#1}\else
 \def\theurl{\previousurls\hskip 0.75em{\rm and}\hskip 0.75em
 \relax#1}\fi}
\long\def\abstract#1\endabstract{\long\def\theabstract{#1}}
\def\primaryclass#1{\def\theprimaryclass{#1}}
\def\secondaryclass#1{\def\thesecondaryclass{#1}}
\def\keywords#1{\def\thekeywords{#1}}
%
%
\let\\\par\let\thetitle\relax\let\theshorttitle\relax
\let\theauthors\relax\let\theshortauthors\relax
\let\theaddress\relax\let\theshortaddress\relax
\let\theemail\relax\let\theurl\relax
\let\theabstract\relax\let\theprimaryclass\relax
\let\thesecondaryclass\relax\let\thekeywords\relax
%
%
%
%
\long\def\maketitlepage{    

\vglue 0.2truein   

%
{\parskip=0pt\leftskip 0pt plus 1fil\def\\{\par\smallskip}{\large
\bf\thetitle}\par\medskip}   

\vglue 0.15truein 

%
{\parskip=0pt\leftskip 0pt plus 1fil\def\\{\par}{\sc\theauthors}
\par\medskip}%
 
\vglue 0.1truein 

%
{\small\parskip=0pt
{\leftskip 0pt plus 1fil\def\\{\par}{\sl\theaddress}\par}
\ifx\theemail\relax\else  
\vglue 5pt \def\\{\stdspace{\rm and}\stdspace} 
\cl{Email:\stdspace\tt\theemail}\fi
\ifx\theurl\relax\else    
\vglue 5pt \def\\{\stdspace{\rm and}\stdspace} 
\cl{URL:\stdspace\tt\theurl}\fi\par}

\vglue 7pt 

{\bf Abstract}

\vglue 5pt

\theabstract

\vglue 7pt 

{\bf AMS Classification numbers}\quad Primary:\quad \theprimaryclass\par

Secondary:\quad \thesecondaryclass

\vglue 5pt 

{\bf Keywords:}\quad \thekeywords

\np  

}    
%
%
\long\def\makeshorttitle{    


%
{\parskip=0pt\leftskip 0pt plus 1fil\def\\{\par\smallskip}{\large
\bf\thetitle}\par\medskip}   

\vglue 0.05truein 

%
{\parskip=0pt\leftskip 0pt plus 1fil\def\\{\par}{\sc\theauthors}
\par\medskip}%
 
\vglue 0.03truein 

%
{\small\parskip=0pt
{\leftskip 0pt plus 1fil\def\\{\par}{\sl\ifx\theshortaddress\relax
\theaddress\else\theshortaddress\fi}\par}
\ifx\theemail\relax\else  
\vglue 5pt \def\\{\stdspace{\rm and}\stdspace} 
\cl{Email:\stdspace\tt\theemail}\fi
\ifx\theurl\relax\else    
\vglue 5pt \def\\{\stdspace{\rm and}\stdspace} 
\cl{URL:\stdspace\tt\theurl}\fi\par}

\vglue 10pt 


{\small\leftskip 25pt\rightskip 25pt{\bf Abstract}\stdspace\theabstract

{\bf AMS Classification}\stdspace\theprimaryclass
\ifx\thesecondaryclass\relax\else; \thesecondaryclass\fi\par
{\bf Keywords}\stdspace \thekeywords\par}
\vglue 7pt
}    
\let\maketitle\makeshorttitle        
%
%

\def\volumenumber#1{\def\thevolumenumber{#1}}
\def\volumeyear#1{\def\thevolumeyear{#1}}
\def\pagenumbers#1#2{\def\startpage{#1}\def\finishpage{#2}}
\def\published#1{\def\publishdate{#1}}
\def\received#1{\def\receiveddate{#1}}
\def\revised#1{\def\reviseddate{#1}}
\let\reviseddate\relax
\volumenumber{X}
\volumeyear{20XX}
\pagenumbers{1}{XXX}
\published{XX Xxxember 20XX}

\long\def\makeagttitle{   
\agt\hfill      
\hbox to 60truept{\vbox to 0pt{\vglue -14truept{\bf [Logo here]}\vss}\hss}
\break
{\small Volume \thevolumenumber\ (\thevolumeyear)
\startpage--\finishpage\nl
Published: \publishdate}

\vglue .2truein

{\parskip=0pt\leftskip 0pt plus 1fil\def\\{\par\smallskip}{\large
\bf\thetitle}\par\medskip}   
\vglue 0.05truein 

%
{\parskip=0pt\leftskip 0pt plus 1fil\def\\{\par}{\sc\theauthors}
\par\medskip}%
 
\vglue 0.03truein 


{\small\leftskip 25truept\rightskip 25truept{\bf Abstract}\stdspace\theabstract

{\bf AMS Classification}\stdspace\theprimaryclass
\ifx\thesecondaryclass\relax\else; \thesecondaryclass\fi\par
{\bf Keywords}\stdspace \thekeywords\par}\vglue 7truept

}   


\def\Addresses{\bigskip
{\small \parskip 0pt \leftskip 0pt \rightskip 0pt plus 1fil \def\\{\par}
\sl\theaddress\par\medskip \rm Email:\stdspace\tt\theemail\par
\ifx\theurl\relax\else\smallskip \rm URL:\stdspace\tt\theurl\par\fi}}

\def\agtart{
\hoffset 14truemm
\voffset 31truemm
\font\phead=cmsl9 scaled 950
\font\pnum=cmbx10 scaled 913
\font\pfoot=cmsl9 scaled 950
\headline{\vbox to 0pt{\vskip -4.5mm\line{\small\phead\ifnum
\count0=\startpage ISSN numbers are printed here
\hfill {\pnum\folio}\else\ifodd\count0\def\\{ }%
\ifx\theshorttitle\relax\thetitle\else\theshorttitle\fi\hfill{\pnum\folio}
\else\def\\{ and }{\pnum\folio}\hfill\ifx\theshortauthors\relax\theauthors
\else\theshortauthors\fi\fi\fi}\vss}}
\footline{\vbox to 0pt{\vglue 0mm\line{\small\pfoot\ifnum\count0=\startpage
Copyright declaration is printed here\hfill\else
\agt, Volume \thevolumenumber\ (\thevolumeyear)\hfill\fi}\vss}}
\let\maketitle\makeagttitle\let\makeshorttitle\makeagttitle}


\def\ifplaintex{\expandafter\ifx\csname documentclass\endcsname\relax}

\def\gtp{{\mathsurround=0pt\it $\cal G\mskip-2mu$eometry \&\ 
$\cal T\!\!$opology $\cal P\!$ublications}}  

\def\recd{{\small Received:\qua\receiveddate\ifx\reviseddate\relax
\else\qquad Revised:\qua\reviseddate\fi\par}} 


\def\lognumber#1{\def\thelognumber{#1}}
\def\volumenumber#1{\def\thevolumenumber{#1}}
\def\volumeyear#1{\def\thevolumeyear{#1}}
\def\papernumber#1{\def\thepapernumber{#1}}
\def\pagenumbers#1#2{\def\startpage{#1}\def\finishpage{#2}}
\def\published#1{\def\publishdate{#1}}

\def\received#1{\def\receiveddate{#1}}
\def\revised#1{\def\reviseddate{#1}}

\def\asciiauthors#1{\def\theasciiauthors{#1}}
\def\asciiaddress#1{\def\theasciiaddress{#1}}

\long\def\asciiabstract#1{\long\def\theasciiabstract{#1}}


\let\\\par\let\thelognumber\relax\let\thevolumenumber\relax
\let\thepapernumber\relax\let\thevolumeyear\relax\let\startpage\relax
\let\finishpage\relax\let\publishdate\relax\let\receiveddate\relax
\let\reviseddate\relax\let\theasciititle\relax
\let\theasciiauthors\relax\let\theasciiaddress\relax
\let\theasciiabstract\relax

\let\theasciiemail\relax


\ifplaintex
\font\logobig=cmssbx10 scaled 3836
\font\logomed=cmssbx10 scaled 2557
\else
\font\logobig=cmssbx10 scaled 4200
\font\logomed=cmssbx10 scaled 2800
\fi

\long\def\makeagttitle{   
\count0=\startpage
\agt\hfill      
\hbox to 45truept{\vbox to 0pt{\vglue -13truept{\logomed A\kern -.37em{\logobig 
T}\kern -.38em G}\vss}\hss}
\break
{\small Volume \thevolumenumber\ (\thevolumeyear)
\startpage--\finishpage\nl
Published: \publishdate}

\vglue .25truein

{\parskip=0pt\leftskip 0pt plus
1fil\def\\{\par\smallskip}{\Large\bf\thetitle}\par\medskip} \vglue
0.05truein

%
{\parskip=0pt\leftskip 0pt plus 1fil\def\\{\par}{\sc\theauthors}
\par\medskip}%
 
\vglue 0.03truein 


{\small\leftskip 25truept\rightskip 25truept{\bf Abstract}\stdspace\theabstract

{\bf AMS Classification}\stdspace\theprimaryclass
\ifx\thesecondaryclass\relax\else; \thesecondaryclass\fi\par
{\bf Keywords}\stdspace \thekeywords\par}\vglue 7truept

}   

\ifplaintex
\hoffset 14truemm
\voffset 31truemm
\font\phead=cmsl9 scaled 950
\font\pnum=cmbx10 scaled 913
\font\pfoot=cmsl9 scaled 950
\headline{\vbox to 0pt{\vskip -4.5mm\line{\small\phead\ifnum
\count0=\startpage ISSN 1472-2739 (on-line) 1472-2747 (printed)
\hfill {\pnum\folio}\else\ifodd\count0\def\\{ }%
\ifx\theshorttitle\relax\thetitle\else\theshorttitle\fi\hfill{\pnum\folio}
\else\def\\{ and }{\pnum\folio}\hfill\ifx\theshortauthors\relax\theauthors
\else\theshortauthors\fi\fi\fi}\vss}}
\footline{\vbox to 0pt{\vglue 0mm\line{\small\pfoot\ifnum\count0=\startpage
\copyright\ \gtp\hfill\else
\agt, Volume \thevolumenumber\ (\thevolumeyear)\hfill\fi}\vss}}
\else
\headsep 23pt
\footskip 35pt
\hoffset -4truemm
\voffset 12.5truemm
\font\lhead=cmsl9 scaled 1050
\font\lnum=cmbx10 
\font\lfoot=cmsl9 scaled 1050
\makeatletter
\def\@oddhead{{\small\lhead\ifnum\count0=\startpage ISSN 1472-2739 
(on-line) 1472-2747 (printed)\hfill {\lnum\number\count0}\else\ifodd\count0
\def\\{ }\ifx\theshorttitle\relax \thetitle \else\theshorttitle\fi\hfill
{\lnum\number\count0}\else\def\\{ and }{\lnum\number\count0}
\hfill\ifx\theshortauthors\relax 
\theauthors\else\theshortauthors\fi\fi\fi}}\def\@evenhead{\@oddhead}
\def\@oddfoot{\small\lfoot\ifnum\count0=\startpage\copyright\ \gtp\hfill\else
\agt, Volume \thevolumenumber\ (\thevolumeyear)\hfill\fi}
\def\@evenfoot{\@oddfoot}
\makeatother
\fi
\let\maketitlepage\makeagttitle
\let\makeshorttitle\maketitlepage
\let\maketitle\maketitlepage


\newwrite\gtoutfile
\long\gdef\makeheadfile{  
{\def\\{, }\def\s{ }
\immediate\openout\gtoutfile head.xxx
\immediate\write\gtoutfile{To: math@arxiv.org}
\immediate\write\gtoutfile{Subject: put OR rep NNNNN:ppppp}
\immediate\write\gtoutfile{--text follows this line--}
\immediate\write\gtoutfile{Proxy-for: \ifx\theasciiauthors\relax
\theauthors\else\theasciiauthors\fi\s<\ifx\theasciiemail\relax\theemail\else\theasciiemail\fi>}
\immediate\write\gtoutfile{\noexpand\\}
\immediate\write\gtoutfile{Authors: \ifx\theasciiauthors\relax
\theauthors\else\theasciiauthors\fi}
{\def\\{ }\immediate\write\gtoutfile{Title: \ifx\theasciititle\relax
\thetitle\else\theasciititle\fi}}
\immediate\write\gtoutfile{Subj-class: GT or SG, GR etc}
\immediate\write\gtoutfile{MSC-class: \theprimaryclass\ifx\thesecondaryclass\relax\else, \thesecondaryclass\fi}
\immediate\write\gtoutfile{Journal-ref: Algebraic and Geometric Topology \thevolumenumber\s
(\thevolumeyear) \startpage-\finishpage}
\immediate\write\gtoutfile{Comments: Published by Algebraic and
Geometric Topology at}
\immediate\write\gtoutfile{\s\s\s  http://www.maths.warwick.ac.uk/agt/AGTVol\thevolumenumber/agt-\thevolumenumber-\thepapernumber.abs.html}
\immediate\write\gtoutfile{\noexpand\\}
\immediate\write\gtoutfile{}
\ifx\theasciiabstract\relax
\immediate\write\gtoutfile{\theabstract}\else
\immediate\write\gtoutfile{\theasciiabstract}\fi
\immediate\write\gtoutfile{}
\immediate\write\gtoutfile{\noexpand\\}
\immediate\write\gtoutfile{}
\immediate\closeout\gtoutfile}}  

\def\maketitlepage{\makeagttitle\makeheadfile}
\let\makeshorttitle\maketitlepage
\let\maketitle\maketitlepage


\def\ifplaintex{\expandafter\ifx\csname documentclass\endcsname\relax}

\def\gtp{{\mathsurround=0pt\it $\cal G\mskip-2mu$eometry \&\ 
$\cal T\!\!$opology $\cal P\!$ublications}}  

\def\recd{{\small Received:\qua\receiveddate\ifx\reviseddate\relax
\else\qquad Revised:\qua\reviseddate\fi\par}} 


\def\lognumber#1{\def\thelognumber{#1}}
\def\volumenumber#1{\def\thevolumenumber{#1}}
\def\volumeyear#1{\def\thevolumeyear{#1}}
\def\papernumber#1{\def\thepapernumber{#1}}
\def\pagenumbers#1#2{\def\startpage{#1}\def\finishpage{#2}}
\def\published#1{\def\publishdate{#1}}

\def\received#1{\def\receiveddate{#1}}
\def\revised#1{\def\reviseddate{#1}}

\def\asciiauthors#1{\def\theasciiauthors{#1}}
\def\asciiaddress#1{\def\theasciiaddress{#1}}

\long\def\asciiabstract#1{\long\def\theasciiabstract{#1}}


\let\\\par\let\thelognumber\relax\let\thevolumenumber\relax
\let\thepapernumber\relax\let\thevolumeyear\relax\let\startpage\relax
\let\finishpage\relax\let\publishdate\relax\let\receiveddate\relax
\let\reviseddate\relax\let\theasciititle\relax
\let\theasciiauthors\relax\let\theasciiaddress\relax
\let\theasciiabstract\relax

\let\theasciiemail\relax


\ifplaintex
\font\logobig=cmssbx10 scaled 3836
\font\logomed=cmssbx10 scaled 2557
\else
\font\logobig=cmssbx10 scaled 4200
\font\logomed=cmssbx10 scaled 2800
\fi

\long\def\makeagttitle{   
\count0=\startpage
\agt\hfill      
\hbox to 45truept{\vbox to 0pt{\vglue -13truept{\logomed A\kern -.37em{\logobig 
T}\kern -.38em G}\vss}\hss}
\break
{\small Volume \thevolumenumber\ (\thevolumeyear)
\startpage--\finishpage\nl
Published: \publishdate}

\vglue .25truein

{\parskip=0pt\leftskip 0pt plus
1fil\def\\{\par\smallskip}{\Large\bf\thetitle}\par\medskip} \vglue
0.05truein

%
{\parskip=0pt\leftskip 0pt plus 1fil\def\\{\par}{\sc\theauthors}
\par\medskip}%
 
\vglue 0.03truein 


{\small\leftskip 25truept\rightskip 25truept{\bf Abstract}\stdspace\theabstract

{\bf AMS Classification}\stdspace\theprimaryclass
\ifx\thesecondaryclass\relax\else; \thesecondaryclass\fi\par
{\bf Keywords}\stdspace \thekeywords\par}\vglue 7truept

}   

\ifplaintex
\hoffset 14truemm
\voffset 31truemm
\font\phead=cmsl9 scaled 950
\font\pnum=cmbx10 scaled 913
\font\pfoot=cmsl9 scaled 950
\headline{\vbox to 0pt{\vskip -4.5mm\line{\small\phead\ifnum
\count0=\startpage ISSN 1472-2739 (on-line) 1472-2747 (printed)
\hfill {\pnum\folio}\else\ifodd\count0\def\\{ }%
\ifx\theshorttitle\relax\thetitle\else\theshorttitle\fi\hfill{\pnum\folio}
\else\def\\{ and }{\pnum\folio}\hfill\ifx\theshortauthors\relax\theauthors
\else\theshortauthors\fi\fi\fi}\vss}}
\footline{\vbox to 0pt{\vglue 0mm\line{\small\pfoot\ifnum\count0=\startpage
\copyright\ \gtp\hfill\else
\agt, Volume \thevolumenumber\ (\thevolumeyear)\hfill\fi}\vss}}
\else
\headsep 23pt
\footskip 35pt
\hoffset -4truemm
\voffset 12.5truemm
\font\lhead=cmsl9 scaled 1050
\font\lnum=cmbx10 
\font\lfoot=cmsl9 scaled 1050
\makeatletter
\def\@oddhead{{\small\lhead\ifnum\count0=\startpage ISSN 1472-2739 
(on-line) 1472-2747 (printed)\hfill {\lnum\number\count0}\else\ifodd\count0
\def\\{ }\ifx\theshorttitle\relax \thetitle \else\theshorttitle\fi\hfill
{\lnum\number\count0}\else\def\\{ and }{\lnum\number\count0}
\hfill\ifx\theshortauthors\relax 
\theauthors\else\theshortauthors\fi\fi\fi}}\def\@evenhead{\@oddhead}
\def\@oddfoot{\small\lfoot\ifnum\count0=\startpage\copyright\ \gtp\hfill\else
\agt, Volume \thevolumenumber\ (\thevolumeyear)\hfill\fi}
\def\@evenfoot{\@oddfoot}
\makeatother
\fi
\let\maketitlepage\makeagttitle
\let\makeshorttitle\maketitlepage
\let\maketitle\maketitlepage


\newwrite\gtoutfile
\long\gdef\makeheadfile{  
{\def\\{, }\def\s{ }
\immediate\openout\gtoutfile head.xxx
\immediate\write\gtoutfile{To: math@arxiv.org}
\immediate\write\gtoutfile{Subject: put OR rep NNNNN:ppppp}
\immediate\write\gtoutfile{--text follows this line--}
\immediate\write\gtoutfile{Proxy-for: \ifx\theasciiauthors\relax
\theauthors\else\theasciiauthors\fi\s<\ifx\theasciiemail\relax\theemail\else\theasciiemail\fi>}
\immediate\write\gtoutfile{\noexpand\\}
\immediate\write\gtoutfile{Authors: \ifx\theasciiauthors\relax
\theauthors\else\theasciiauthors\fi}
{\def\\{ }\immediate\write\gtoutfile{Title: \ifx\theasciititle\relax
\thetitle\else\theasciititle\fi}}
\immediate\write\gtoutfile{Subj-class: GT or SG, GR etc}
\immediate\write\gtoutfile{MSC-class: \theprimaryclass\ifx\thesecondaryclass\relax\else, \thesecondaryclass\fi}
\immediate\write\gtoutfile{Journal-ref: Algebraic and Geometric Topology \thevolumenumber\s
(\thevolumeyear) \startpage-\finishpage}
\immediate\write\gtoutfile{Comments: Published by Algebraic and
Geometric Topology at}
\immediate\write\gtoutfile{\s\s\s  http://www.maths.warwick.ac.uk/agt/AGTVol\thevolumenumber/agt-\thevolumenumber-\thepapernumber.abs.html}
\immediate\write\gtoutfile{\noexpand\\}
\immediate\write\gtoutfile{}
\ifx\theasciiabstract\relax
\immediate\write\gtoutfile{\theabstract}\else
\immediate\write\gtoutfile{\theasciiabstract}\fi
\immediate\write\gtoutfile{}
\immediate\write\gtoutfile{\noexpand\\}
\immediate\write\gtoutfile{}
\immediate\closeout\gtoutfile}}  

\def\maketitlepage{\makeagttitle\makeheadfile}
\let\makeshorttitle\maketitlepage
\let\maketitle\maketitlepage

\lognumber{2}
\volumenumber{1}
\volumeyear{2001}
\papernumber{2}
\published{15 November 2000}
\pagenumbers{31}{37}
\received{8 August 2000}\revised{9 October 2000}

\def\Cal#1{{\cal #1}}

\reflist

\key {E}
{\bf M Eudave-Mu\~noz},
{\it Meridional essential surfaces for tunnel number one knots},
To appear in Bol. Soc. Mat. Mex. (third series) 6 (2000).

\key {FF}
{\bf B Freedman}, {\bf M\,H Freedman},
{\it Kneser--Haken finiteness for bounded  3-mani-\break folds, locally free
groups, and cyclic covers},
Topology {37} (1998) 133--147.

\key {H}
{\bf J Hempel},
{\it 3-manifolds},
Annals of Math. Studies {86}, Princeton University Press (1976).

\key {Ho}
{\bf  H Howards},
{\it Arbitrarily many incompressible surfaces in three-manifolds with genus two
boundary components},
preprint.

\key {J}
{\bf W Jaco},
{\it  Lectures on Three-Manifold Topology}, 
Regional Conference Series in Mathematics 43, 
Amer. Math. Soc. (1981).

\key {K} 
{\bf T Kobayashi},
{\it Structures of full Haken manifolds},
Osaka J. Math {24}  (1987) 173--215.

\key {M}
{\bf K Morimoto},
{\it Tunnel number, connected sum and meridional essential surfaces},
Topology {39} (2000) 469--485.

\key {SS} 
{\bf M Scharlemann}, {\bf J Schultens},
{\it Comparing Heegaard and JSJ structures of orientable
3-manifolds}, to appear in Trans. Amer. Math. Soc.

\key {ST}
{\bf M Scharlemann}, {\bf A Thompson},
{\it Thin position for 3-manifolds},
AMS Contemporary Mathematics {164} (1994) 231--238.

\key {Sc} 
{\bf J Schultens},
{\it Additivity of tunnel number for small knots},
preprint.

\key {Sh}
{\bf W Sherman},
{\it On the Kneser--Haken finiteness theorem: a sharpness result},
Ph.D. Thesis, University of California, Los Angeles, 1992.

\endreflist

\title {A universal bound for surfaces in 3-manifolds\\with a given Heegaard
genus}
\shorttitle{A universal bound for surfaces in 3-manifolds}

\authors {Mario Eudave-Mu\~noz\\Jeremy Shor}
\asciiauthors {Mario Eudave-Munoz and Jeremy Shor}

\address {Instituto de Matem\'aticas, UNAM, Circuito Exterior, Ciudad
Universitaria\\04510 M\'exico D.F., MEXICO\\\smallskip\\and
\\\smallskip\\250 E. Houston St.~no.~3J, New York City, NY 10002, U.S.A.}
\asciiaddress {Instituto de Matematicas, UNAM, Circuito Exterior, Ciudad
Universitaria\\04510 Mexico D.F., MEXICO\\250 E. Houston St. 
no. 3J, New York City, NY 10002, U.S.A.}

\email{mario@matem.unam.mx}

\abstract It is shown that for given positive integers $g$ and $b$,
there is a number $\Cal C (g,b)$, such that any orientable compact
irreducible 3-manifold of Heegaard genus $g$ has at most 
$\Cal C (g,b)$ disjoint, nonparallel incompressible surfaces with first
Betti number $b_1 < b$. 
\endabstract
\asciiabstract{It is shown that for given positive integers g and b,
there is a number C(g,b), such that any orientable compact irreducible
3-manifold of Heegaard genus $g$ has at most C(g,b) disjoint,
nonparallel incompressible surfaces with first Betti number b_1 < b.}

\primaryclass{57N10}
\secondaryclass{57M25}
\keywords{Incompressible surface, Haken manifold.} 

\makeshorttitle

\section{Introduction}

Let $M$ be a compact, orientable, irreducible 3-manifold, possibly
with boundary.  As usual, by a surface we mean a connected, compact
2-manifold; if it is contained in a 3-manifold, then it is assumed to
be properly embedded.  There are several results which bound the
number of disjoint incompressible surfaces in $M$.

The Haken--Kneser finiteness theorem says that given $M$, there exist
an integer $c(M)$, such that any collection of pairwise disjoint,
non-parallel, closed, incompressible surfaces in $M$ has at most
$c(M)$ components.  This can be extended to surfaces with boundary, if
the surfaces are taken to be incompressible and
$\partial$-incompressible (for a proof see [\ref {H}] and [\ref {J}]).

However, this theorem is false if the surfaces are not
$\partial$-incompressible. Howards [\ref {Ho}] has shown that if a
3-manifold $M$ has a boundary component of genus 2 or greater, then
there are in $M$ arbitrary large collections of disjoint,
non-parallel, incompressible surfaces. This was proved previously by
Sherman [\ref {Sh}] for $F \times I$, where $F$ is a closed surface of
sufficiently high genus. But in these constructions, when there are
large numbers of surfaces, the Betti numbers of the surfaces are also
large. If we bound the Betti numbers, then the number of surfaces will
be bounded by in the following theorem, proved by Freedman and
Freedman [\ref {FF}].

\proclaim {Theorem 1}{\rm(Freedman and Freedman [\ref {FF}])}\qua 
Let $M$ be a compact 3-manifold with boundary and $b$ an integer
greater than zero. There is a constant $c(M,b)$ so that if $F_1,\dots,
F_k$, $k>c$ is a collection of disjoint, incompressible surfaces such
that all Betti numbers $b_1(F_i)<b$, $1\leq i \leq k$, and no $F_i$,
$1\leq i \leq k$, is a boundary parallel annulus or a boundary
parallel disk, then at least two members $P_i$ and $P_j$ are parallel.
\endproclaim

However in all of these theorems, the bound depends on a given
manifold, i.e., the bound is not universal.  We intend to look for
some kind of universal bound.  For example, if we restrict to the
class of closed, irreducible 3-manifolds, then as said above, each
manifold has a bounded number of incompressible surfaces, but clearly
there is no universal bound, for it is easy to construct closed
manifolds with a large number of incompressible surfaces, even with
surfaces of a given genus.  If we bound the Heegaard genus of the
3-manifolds, then there is no such bound, for in [\ref {E}] a
construction is given, so that for each integer $N$, we can find a
closed, irreducible 3-manifold with Heegaard genus 2 and having more
than $N$ disjoint, non-parallel, incompressible surfaces.  But in
those examples, when the number of surfaces is high, the genus of the
surfaces is also high. If we bound also the genus of the surfaces,
say, if we restrict to tori, then it follows from [\ref {K}] and [\ref
{SS}], that there is an integer $N(g)$, depending only on $g$, such
that any closed 3-manifold with Heegaard genus $g$, has at most $N(g)$
disjoint, non-parallel, incompressible tori.  All this suggests that a
result such as the following could exist.  This is an extension of
Theorem 1 to all manifolds of a given Heegaard genus.

\proclaim {Theorem 2}{\rm(Main theorem)}\qua 
Consider the class $\Cal M^3(g)$ of all compact, orientable and
irreducible 3-manifolds with Heegaard genus $g$. Let $b$ be an integer
greater than zero. Then there is a constant $\Cal C(g,b)$, depending
only on $g$ and $b$, such that for any $M$ in $\Cal M^3(g)$, if
$P_1,\dots, P_k$, $k> \Cal C(g,b)$ is a collection of disjoint,
incompressible surfaces in $M$, such that all Betti numbers
$b_1(P_i)<b$, $1\leq i \leq k$, and no $P_i$, $1\leq i
\leq k$, is a boundary parallel annulus or a boundary parallel disk, then at
least two members $P_i$ and $P_j$ are parallel.
\endproclaim

Here is an outline of the proof.  We consider a 3-manifold $M$ with a
Heegaard splitting of genus $g$. If the splitting is not strongly
irreducible, we untelescope it [\ref {ST}], decomposing $M$ into a
union of submanifolds, glued along incompressible surfaces, and such
that each submanifold has a strongly irreducible Heegaard
splitting. Given a collection of surfaces in $M$, these can be
isotoped to intersect the surfaces of the decomposition of $M$
essentially. Then we can apply Theorem 1 to each of the pieces, which
is a compression body of genus at most $g$, and by a counting
argument, similar to the one in the Haken--Kneser Theorem which counts
bad pieces in a tetrahedra, but now counting bad pieces in a
compression body, we get the required bound.

For closed manifolds we deduce.

\proclaim {Corollary 3} Consider the class $\Cal M^3_c(g)$ of all closed,
orientable and irreducible 3-manifolds with Heegaard genus $g$. Let
$h$ be an integer greater than zero. Then there is a constant $C(g,h)$
depending only on $g$ and $h$, such that for any $M$ in $\Cal
M^3_c(g)$, if $P_1,\dots, P_k$, $k>C(g,h)$ is a collection of
disjoint, incompressible surfaces in $M$, all of genus less than $h$,
then at least two members $P_i$ and $P_j$ are parallel.
\endproclaim

\rk{Acknowledgements}This work was done while the first author was 
visiting the University of California at Santa Barbara, and he was
partially supported by a NSF grant and by a sabbatical fellowship from
CONACYT (Mexico). The first author is grateful to Eric Sedgwick by
suggesting him that something like Corollary 3 could be true.

\section {Proof and estimates}

A compression body is a 3-manifold $W$ obtained from a connected
closed orientable surface $S$ by attaching 2-handles to $S\times \{
0\}\subset S\times I$ and capping off any resulting 2-sphere boundary
components. We denote $S\times \{1\}$ by $\partial_+ W$ and $\partial
W-\partial_+W$ by $\partial_- W$. Dually, a compression body is a
connected, orientable 3-manifold obtained from a (not necessarily
connected) closed orientable surface $\partial_-W\times I$ by
attaching 1-handles. Define the index of $W$ by $J(W)=\chi
(\partial_-W)-\chi (\partial_+W) \geq 0$. The genus of a compression
body $W$ is defined to be the genus of the surface $\partial_+ W$.

A Heegaard splitting of a 3-manifold $M$ is a decomposition $M=V\cup_S
W$, where $V$, $W$ are compression bodies such that $V\cap
W=\partial_+V=\partial_+W=S$. We call $S$ the splitting surface or
Heegaard surface. We say that the splitting is of genus $g$ if the
surface $S$ is of genus $g$. We say that a 3-manifold has Heegaard
genus $g$ if $M$ has a splitting of genus $g$, and any other splitting
has genus $\geq g$.

A Heegaard splitting is reducible if there are essential disks $D_1$
and $D_2$, in $V$ and $W$ respectively, such that $\partial
D_1=\partial D_2$.  A splitting is irreducible if it is not reducible.
A Heegaard splitting is weakly reducible if there are essential disks
$D_1$ and $D_2$, in $V$ and $W$ respectively, such that $\partial D_1
\cap \partial D_2 = \emptyset$.  A Heegaard splitting is strongly
irreducible if it is not weakly reducible.

A generalized Heegaard splitting of a compact orientable 3-manifold $M$ is a
structure $M=(V_1\cup_{S_1} W_1)\cup_{F_1}(V_2\cup_{S_2} W_2)\cup_{F_2}\dots
\cup_{F_{m-1}} (V_m\cup_{S_m} W_m)$. Each of the $V_i$ and $W_i$ is a union
of compression bodies, $\partial_+ V_i=S_i=\partial_+W_i$, (i.e.,
$V_i\cup_{S_i} W_i$ is a union of Heegaard splittings of a submanifold
of M) and $\partial_- W_i=F_i=\partial_-V_{i+1}$. We say that a
generalized Heegaard splitting is strongly irreducible if each
Heegaard splitting $V_i\cup_{S_i} W_i$ is strongly irreducible and
each $F_i$ is incompressible in $M$. We will denote $\cup_i F_i$ by
$\Cal F$ and $\cup_i S_i$ by $\Cal S$. The surfaces in $\Cal F$ are
called the thin levels and the surfaces in $\Cal S$ the thick levels.

\proclaim {Theorem A} Let $M$ be a compact irreducible 3-manifold.
Suppose the Heegaard genus of $M$ is $g$. Then $M$ has a strongly
irreducible generalized Heegaard splitting $M=(V_1\cup_{S_1}
W_1)\cup_{F_1}(V_2\cup_{S_2} W_2)\cup_{F_2}\dots
\cup_{F_{m-1}} (V_m\cup_{S_m} W_m)$,
such that $\sum J(V_i)= \sum J(W_i) = 2g-2$.
\endproclaim

A proof of this Theorem can be found in [\ref {ST}] and [\ref
{SS}]. The following Theorem is proved in [\ref {Sc}]. A different
proof is given in [\ref {M}].

\proclaim {Theorem B} Let $P$ be a properly embedded, not necessarily 
connected, incompressible surface in an irreducible 3-manifold $M$,
and let $M=(V_1\cup_{S_1} W_1)\cup_{F_1}(V_2\cup_{S_2}
W_2)\cup_{F_2}\dots\cup_{F_{m-1}} (V_m\cup_{S_m} W_m)$ be a strongly
irreducible generalized Heegaard splitting of $M$. Then $\Cal F\cup
\Cal S$ can be isotoped to intersect $P$ only in curves that are
essential in both $P$ and $\Cal F\cup \Cal S$.
\endproclaim

\proof {Proof of Theorem 2}

Let $M$ be an orientable, irreducible, compact 3-manifold of Heegaard
genus $g$, that is, there are compression bodies $V$ and $W$ in $M$,
with $S=\partial_+V = \partial_+W = V\cap W$, $M=V\cup_S W$, and
genus$(S)=g$.

By Theorem A, $M$ possesses a strongly irreducible Heegaard splitting
$M=(V_1\cup_{S_1} W_1)\cup_{F_1}(V_2\cup_{S_2} W_2)\cup_{F_2}\dots
\cup_{F_{m-1}} (V_m\cup_{S_m} W_m)$, such that $\sum J(V_i)= \sum J(W_i) 
= 2g-2$.  For a fixed $i$, $V_i\cup_{S_i} W_i$ is a union of strongly
irreducible Heegaard splittings, only one component of it is not a
product, this is called the active component. If $V_{ij}$ is a
component of $V_i$, then $J(V_{ij})=0$, only if $V_{ij}$ is a product,
and $J(V_{ij})\geq 2$ otherwise. It follows that $m\leq g-1$, and that
there are at most $g-1$ active components.

Let $P_1, P_2, \dots, P_n$ be a collection of disjoint, properly
embedded, incompressible surfaces in $M$. Suppose that $b_1 (P_i) <
b$, for all $i$.  By Theorem B, $\Cal F \cup S$ can be isotoped to
intersect $P_1, P_2, \dots, P_n$ only in curves which are essential in
both $P_1, P_2, \dots, P_n$ and $\Cal F \cup S$.  If some component of
$V_i\cup_{S_i} W_i$ is a product, say, it is of the form $S\times
[-1,1]$, then by taking the product sufficiently small, we can suppose
that $(S\times [-1,1])\cap P_j=((S\times \{0\})\cap P_j )\times
[-1,1]$ for any surface $P_j$. So we can disregard the product
$S\times [-1,1]$, and consider only the surface $S\times\{1\}$. Doing
this with all product components, we can assume, without loss of
generality, that the decomposition of $M$ consist only of the active
components. So $M$ is decomposed into at most $2g-2$ compression
bodies.

Consider $P_j\cap V_i$; this is a collection of surfaces. Each of the
components has $b_1 < b$, for $P_j$ intersect $\Cal S$ and $\Cal F$ in
simple closed curves which are essential in $P_j$. We can assume that
no component of $P_j\cap V_i$ is a boundary parallel annulus in $V_i$,
for if there is one, then by an isotopy of $P_i$, it can be pushed
into another compression body, reducing the number of intersections
between $P_j$ and $\Cal F \cup \Cal S$.

For each $V_i$ (and $W_i$), by Theorem 1, there is a constant
$c(V_i,b)$ depending only on $V_i$ and $b$, such that any collection
of incompressible surfaces in $V_i$, with $b_1 < b$, no two of them
parallel, has at most $c(V_i,b)$ components.

Let $W$ be a compression body. Let $F_1,\dots,F_k$ a collection of
incompressible surfaces in $W$. A component of $W-(F_1\cup\dots \cup
F_k)$ is good, if it is of the form $F_i\times I$, $I=[0,1]$, where
$F\times \{0\}$ and $F\times \{1\}$ are surfaces in the collection,
and $F\times \partial I \subset \partial W$.  A component of
$W-(F_1\cup\dots \cup F_k)$ is bad, if it is not good.  If all the
$F_i's$ satisfy that $b_1(F_i) < b$, then it follows from Theorem 1,
that for any such collection, $W-(F_1\cup\dots \cup F_k)$ has at most
$c(W,b) +1$ bad components.

Let $C=max \{ c(V,b); V$ a compression body of genus $\leq g  \}$. And let
$$\Cal C(g,b) = 4g + (2g-2)(C+1)$$
Note that $\Cal C(g,b)$ depends only on $g$ and $b$.  For the
collection of surfaces $P_1,\dots,$ $P_n$, suppose further that $n >
\Cal C(g,b)$.  As $M$ has Heegaard genus $g$, then $b_1(M)\leq
2g$. This implies that $M-(P_1\cup P_2\cup \dots \cup P_n)$ has at
least $\Cal C(g,b)-2g +1$ components.  In each $V_i$ and $W_i$ there
are at most $C+1$ bad pieces, so there are in total at most
$(2g-2)(C+1)$ bad pieces. Then at least $2g+1$ components of
$M-(P_1\cup P_2\cup \dots \cup P_n)$ are made of good pieces. A
component made of good pieces is a $I$-bundle over a surface, but
again $b_1(M,Z_2) \leq 2g$, so at most $2g$ of them can be $I$-bundles
over a nonorientable surface. Then at least a component of $M-(P_1\cup
P_2\cup \dots
\cup P_n)$ is a product, which implies that some $P_i$ and $P_j$ are parallel.
This completes the proof of Theorem 2.
\endproof

\rk{Estimates}

A rough estimate for $\Cal C(g,b)$ can be obtained from [\ref {FF}]. From
that paper it follows that if $V$ is a compression body of genus $g$,
which is not a product, then  a bound $c(V,b)$ is given by 
$$c(V,b)={14\over 3}c_03^{2b-2} +5\vert \chi(V)\vert +T + b_1(V,Z_2)$$
where $c_0$ is the constant given by the Kneser-Haken Theorem, 
and $T$ is the number of tori boundary components. An incompressible and
$\partial$-incompressible surface in a compression body $V$ is either an
essential disk, or an annulus with one boundary component on $\partial_+V$
and the other in $\partial_-V$. 
It is not difficult to see that if $V$ has genus
$g$ then $c_0(v)=3g-3$. Also $T\leq g$, $\vert \chi(V)\vert\leq 2g-3$, and
$b_1(V,Z_2)\leq 2g$.

It follows that
$$c(V,b)\leq 14(g-1)3^{2b-2} +5(2g-3) +g + 2g=14(g-1)3^{2b-2} +13g-15.$$

Define $C=14(g-1)3^{2b-2} +13g-15$. So $c(W,b)\leq C$ for any
compression body $W$ of genus $\leq g$. Then $C$ can be taken to be the
constant
defined in the proof of Theorem 2. From this follows that the constant
$\Cal C(g,b)$ can be expressed as
$$\Cal C(g,b)=p_1(g) + p_2(g) 3^{2b-2}$$
where $p_1(g)$ and $p_2(g)$ are quadratic polynomials in $g$.

\rk{Remarks} Theorem 2 can be extended to arbitrary orientable compact
3-manifolds. This follows from the existence of a decomposition on prime
3-manifolds, the additivity of the Heegaard genus, and from the fact that given
any collection of incompressible surfaces in a 3-manifold $M$, then there is a
complete collection of decomposing spheres in $M$ which are disjoint from the
given surfaces. 

An extension of Theorem 2 for nonorientable manifolds would follow
from a version of Theorems A and B for nonorientable 3-manifolds.

\references

\Addresses
\recd

\bye